\documentclass[12 pt]{article}

\usepackage{latexsym}
\usepackage{amssymb}
\usepackage{amsmath}
\usepackage{amsfonts}
\usepackage{mathrsfs}
\usepackage[all]{xy}
\usepackage{texdraw}
\usepackage{graphicx}
\usepackage{psfrag}
\usepackage{makeidx}

\newtheorem{Theory}{Theory}[section] 
\newtheorem{theorem}[Theory]{Theorem}
\newtheorem{lemma}[Theory]{Lemma}
\newtheorem{technicalLemma}[Theory]{Technical Lemma}
\newtheorem{corollary}[Theory]{Corollary}
\newtheorem{proposition}[Theory]{Proposition}
\newtheorem{fact}{Fact}  
\newtheorem{remark}[Theory]{Remark} 
\newtheorem{question}{Question} 
\newtheorem{conjecture}[question]{Conjecture}

\newtheorem{Ntn}{Description} 
\newtheorem{Dn}[Ntn]{Definition}

\newcommand{\be}{\begin{enumerate}}
\newcommand{\ee}{\end{enumerate}}
\newcommand{\bq}{\begin{question}}
\newcommand{\eq}{\end{question}}
\newcommand{\bcj}{\begin{conjecture}}
\newcommand{\ecj}{\end{conjecture}}
\newcommand{\bc}{\begin{corollary}}
\newcommand{\ec}{\end{corollary}}
\newcommand{\bl}{\begin{lemma}}
\newcommand{\el}{\end{lemma}}
\newcommand{\btl}{\begin{technicalLemma}}
\newcommand{\etl}{\end{technicalLemma}}
\newcommand{\bt}{\begin{theorem}}
\newcommand{\et}{\end{theorem}}
\newcommand{\bp}{\begin{proposition}}
\newcommand{\ep}{\end{proposition}}
\newcommand{\bft}{\begin{fact}}
\newcommand{\eft}{\end{fact}}
\newcommand{\brk}{\begin{remark}}
\newcommand{\erk}{\end{remark}}
\newcommand{\bd}{\begin{Dn}}
\newcommand{\ed}{\end{Dn}}

\newcommand{\ga}{\alpha}
\newcommand{\gb}{\beta}

\newcommand{\ploi}{PL_o(I) }

\newcommand{\Z}{ \mathbf Z }
 
\newcommand{\N}{\mathbf N }
 
\newcommand{\R}{\mathbf R }

\newcommand{\mm}{\mathscr{M}}
\newcommand{\ws}{\mathscr{R}}

\newcommand{\cg}{\mathscr{G}}
\newcommand{\bsc}{\mathcal{B}}
\newcommand{\wc}{\mathcal{W}}
\renewcommand{\sc}{\mathcal{S}}
\newcommand{\pc}{\mathcal{P}}

\newcommand{\rtr}{ \R\times\R }

\newcommand{\bSeq}[4]{\left\{#1_{#2}\right\}_{#2=#3}^{#4}}
\newcommand{\pSeq}[4]{(#1_{#2})_{#2=#3}^{#4}}

\newcommand{\pow}[1]{\mathscr{P}(#1)}

\newcommand{\supp}{\operatorname{Supp}}

\author{Collin Bleak}

\usepackage{epsfig}

\begin{document}
\centerline{\LARGE An Algebraic Classification of Some Solvable
Groups} 
\centerline{\LARGE of Homeomorphisms}

\vspace{.3 in}

{\centerline {\large Abstract}} We produce two separate algebraic
descriptions of the isomorphism classes of the solvable subgroups of
the group $\ploi$ of piecewise-linear orientation-preserving
homeomorphisms of the unit interval under the operation of
composition, and also of the generalized R. Thompson groups $F_n$.
The first description is as a set of isomorphism classes of groups
which is closed under three algebraic operations, and the second is as
the set of isomorphism classes of subgroups of a countable collection of
easily described groups.  We show the two descriptions are equivalent.

\vspace{.1 in}

\tableofcontents

\section{Introduction}
We use the main geometric result of \cite{bpgsc} to produce two
algebraic descriptions of the solvable subgroups of $\ploi$, the group
of orientation-preserving piecewise-linear homeomorphisms of the unit
interval with finitely many breaks in slope under the operation of
composition.  Our results apply equally to any group in the family
$\left\{F_n\right\}$ of generalized Thompson groups, each of which has
a definition as a particular subgroup of $\ploi$.    (The groups
$F_n$ were introduced by Brown in \cite{BrownFinite}, where they were
denoted $F_{n,\infty}$. These groups were later extensively studied by
Stein in \cite{SteinPLGroups}, by Brin and Guzm\'an in
\cite{BGAutomorphisms} and by Burillo, Cleary, and Stein in
\cite{BCS}.)

The first description that we find for the isomorphism classes of the
solvable subgroups of $\ploi$ is as the smallest non-empty class $\ws$
of isomorphism classes of groups which is closed under three
operations, \be
\item taking subgroups,
\item forming standard restricted wreath products with $\Z$ (i.e., $K\mapsto
K\wr\Z$), and 
\item forming bounded direct sums (countable direct sums of groups in
the class with a universal bound on their derived length).  
\ee 

(Note that we give these operations as operations on groups, not on
isomorphism classes of groups.  We will generally not track the
distinction between a group and its isomorphism class in discussion,
in order to simplify our language.  An example of
this is that we may write that $\mm\subset\ws$, where $\mm$ is a set
of groups, all of whose isomorphism classes are elements of $\ws$.  If
$H$ is a group, we may also write that $H\in \ws$ when in fact the
isomorphism class of $H$ is an element of $\ws$.  This practice should
cause the reader no confusion, and will not effect our results.  The
statements of results in this section will not follow this convention
and are formally correct.)

\bt
\label{solveClassification}

$H$ is the isomorphism class of a solvable subgroup of $\ploi$ if and
only if $H\in\ws$.

\et

For each natural number $n$, define a group $G_n$ according to the
process below.  First, define $G_0 = 1$, the trivial
group.  Now inductively define
\[
G_n = \bigoplus_{k\in\Z}(G_{n-1}\wr\Z),
\]
for all positive integers $n$.  Let $\mm=
\left\{G_n\,|\,n\in\N\right\}$.  Our second description is as
indicated by the theorem below.

\bt
\label{RandM}

$H\in\ws$ if and only if $H$ is the isomorphism class of a subgroup of
a group in $\mm$.  \et

In particular, $G$ is a solvable subgroup of $\ploi$ if and only if
$G$ is isomorphic to a subgroup of a group in $\mm$.  The proof of the
above theorem will depend on the following lemma.

\bl
\label{solveInM}
If $G$ is a solvable subgroup of $\ploi$ with derived
length $n$, then $G$ is isomorphic to a subgroup of $G_n$.  
\el

It is immediate from construction that the groups in $\mm$ are all
countable, so the last lemma also has the following corollary.

\bc
If $H$ is a solvable subgroup of $\ploi$, then $H$ is countable.
\ec

Finally, each group $F_n$ contains isomorphic copies of each of the
groups in $\mm$, so that we have the following consequence to the
above results.

\bc 

Let $n$ be any integer with $n\geq 2$, and let $H$ be an isomorphism
class of groups.  $H\in\ws$ if and only if $H$ is the isomorphism
class of a solvable subgroup of $F_n$.

\ec

In \cite{NavasSolv}, \cite{NavasPL}, Navas uses dynamics to
analyze the solvable subgroups of a group of homeomorphisms.  His work
there is focussed on the group $Diff^2_+(\R)$.  In \cite{NavasPL}, he
applies his results to show that a finitely generated solvable
subgroup of $\ploi$ with connected support is isomorphic to a
semi-direct product of a group $H$ with the integers $\Z$, where $H$
is a group in a particular class of groups.  This result is contained
in our investigations below, although his techniques are quite
different from our own.

In the next paragraph, I give a partial list of other related work.
This list spans from work directly related to the group $\ploi$ to less
directly related work in the theory of groups acting on one
dimensional manifolds and the theory of the Godbillon-Vey invariant.

We begin by mentioning the work of Brin and Squier, and later of Brin,
on piecewise linear groups of homeomorphisms on the line and the unit
interval, \cite{BSPLR, picric, BrinU, BrinEG}.  We will refer to this
work throughout the paper.  Also directly impacting the theory of
$\ploi$ are the works of Tsuboi, Minakawa, and Oikonomides.  In
\cite{SmallCommutators}, Tsuboi investigates aspects of $\ploi$ while
pursuing his analysis of the Godbillon-Vey invariant, and in
\cite{TTransvections}, he gives another description of the
Higman-Thompson group $T$ via restricted actions of $SL(2;\Z)$ on the
unit circle.  Minakawa in \cite{ExoticCirc1, ExoticCirc2} provides an
invariant that classifies centralizers in $\ploi$, using different
techniques than those of by Brin and Squier in \cite{picric}.
Oikonomides in \cite{Oikonomides} also analyzes the Godbillon-Vey
invariant, and her dissertation contains Minakawa's characterization
of the centralizers of elements of $\ploi$ as well.  In another
direction, Burslem and Wilkinson in \cite{BurslemWilkinson} classify
the solvable subgroups of the group of analytic diffeomorphisms of
$S^1$.  In \cite{FF1, FF2, FF3}, Farb and Franks analyze groups of
homeomorphisms of one-manifolds in analogy with the theory of Lie
groups and their discrete subgroups.  Finally, the reader is referred
to the survey \cite {GhysCircleActions} by Ghys for more information
about groups acting on the circle.

The author would like to thank Matt Brin for asking for an algebraic
interpretation of the geometric results of \cite{bpgsc}, which lead to
all of the results in this paper.  The author would also like to thank
\'Etiene Ghys, Andres Navas, and Takashi Tsuboi for mentioning various
related works.

The results in this paper are contained in the author's dissertation
written at Binghamton University, although some of the proofs here are
new.

\subsection{Key examples}
In this section, we will mention some key examples.  All of these
examples can be realized in $\ploi$, and in any generalized Thompson
group $F_n$. (We will only specifically realize these groups in
R. Thompson's group $F=F_2$ and in $\ploi$; realizations of our
examples in the other groups $F_n$ are easy to find by making obvious
changes in the constructions below.)

Our examples rely on an understanding of the standard restricted
wreath product.  The reader is referred to the paper by P. M. Neumann
\cite{NeumannW} or the book by J. D. P. Meldrum \cite{Meldrum} for
detailed discussions of these products.  For instance, the reader who
has read the first section of the first chapter of \cite{Meldrum} will
be well prepared for the discussion in this paper.  Nonetheless, a
working definition of a standard restricted wreath product of groups
is given below.

Let $A$ and $T$ be groups.  The restricted wreath product $A\wr T$ can
be thought of as the semidirect product $B\rtimes T$ where $B=
\bigoplus_{t\in T}A$ and where the copy of $T$ in the semi-direct
product acts on $B$ by right multiplication on the index in the direct
sum.  In this context $T$ is called the top group in the product, $A$
is called the bottom group, and $B$ is called the base group. We may
identify the groups $B$ and $T$ with their natural images in the semidirect
product definition without comment.

Taking a restricted wreath product with $\Z$ is something that is easy
to realize in $\ploi$, which is why this activity plays a key role in
the definition of the groups in the class $\mm$, and in the class
$\ws$.  In fact, it is so natural that another collection of groups
becomes relevant. Define $W_0 = 1$, the trivial group, and for all
$i\in\N$, define $W_i = W_{i-1}\wr\Z$, so that
\[
W_i =
(\ldots(((\Z\wr\Z)\wr\Z)\wr \Z)\ldots)\wr\Z,
\]
where there are $i$ appearances of $\Z$ on the right.  Theorem
\ref{RandM} can be rephrased in terms of the $W_i$ instead of the
$G_i$.  However, the $W_i$ lack the corresponding first and third
properties of the $G_i$ stated below, and these are serious
deficiencies in the class, from a computational point of view.

We leave the first property below to the reader. We will
prove the latter two later.

\brk
\label{rkGn}

\be
\item Given any non-negative $i\in\Z$, $G_i\cong \bigoplus_{j\in\Z}
G_i$\qquad(note: the subscript is not the sum index).
\item Given non-negative $n\in\Z$, $G_n$ has derived length $n$.
\item If $H$ is a subgroup of $G_k$ for some non-negative $k\in\Z$, and $H$ has
derived length $n$, then $H$ is isomorphic to a subgroup of $G_n$.
\ee
\erk

\subsection{Geometry}

In this section, we will realize the $W_i$ in $F\leq\ploi$, and use
these realizations to motivate some geometric definitions which will
enable us to state the main geometric result of \cite{bpgsc}.  Here, we
are thinking of $F$ as the realization of Thompson's group in $\ploi$
which consists of all the elements of $\ploi$ which have all slopes
powers of two, and which have all breakpoints occuring at the dyadic
rationals $\Z[\frac{1}{2}]$.  See Cannon, Floyd, and Parry
\cite{CFP} for an introduction to the remarkable group $F$.
\subsubsection{Realizing the $W_i$}
Consider the two elements $\alpha_1$,
$\alpha_2\in\ploi$ defined below:
\[
\begin{array}{llcrr}
x\alpha_1=&\left\{
\begin{array}{lr}
2x&0\leq x<\frac{1}{4},
\\
\\ 
x +\frac{1}{4} & \frac{1}{4}\leq x< \frac{1}{2},
\\
\\ 
\frac{1}{2}x+ \frac{1}{2} &\frac{1}{2}\leq x\leq 1,
\end{array}
\right.&&
x\alpha_2 =&\left\{
\begin{array}{lr}
x&0\leq x<\frac{1}{4},
\\
\\
2x-\frac{1}{4}&\frac{1}{4}\leq x<\frac{5}{16},
\\
\\
x +\frac{1}{16}&\frac{5}{16}\leq x\leq \frac{3}{8},
\\
\\
\frac{1}{2}x +\frac{1}{4}& \frac{3}{8}\leq x< \frac{1}{2},
\\
\\
x &\frac{1}{2}\leq x\leq 1.
\end{array}
\right.
\end{array}
\]
Here are the graphs (superimposed) of these functions:

\begin{center}
\psfrag{a1}[c]{$\alpha_1$}
\psfrag{a2}[c]{$\alpha_2$}
\includegraphics[height=340pt,width=340 pt]{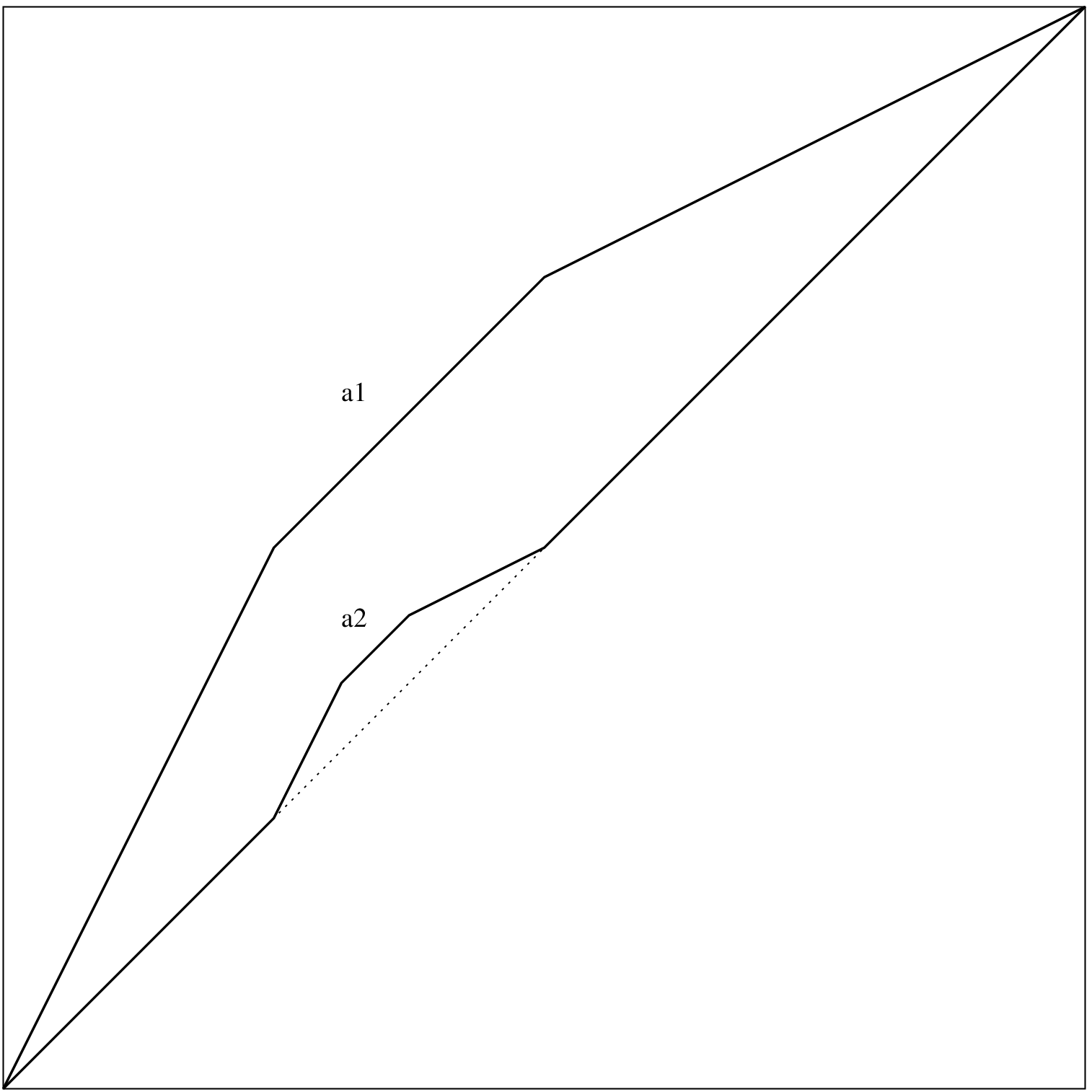}
\end{center}

Either element alone generates a group isomorphic to $\Z\cong W_1$ in
$\ploi$, but the ``action'' of $\alpha_2$ occurs in a single
fundamental domain of $\alpha_1$; that is, $\frac{1}{4}\alpha_1 =
\frac{1}{2}$, but $\alpha_2$ is the identity off of the interval
$[\frac{1}{4},\frac{1}{2}]$.  In particular, $\alpha_2^{\alpha_1} =
\alpha_1^{-1}\alpha_2\alpha_1$ has support
$(\frac{1}{2},\frac{3}{4})$, which is disjoint from the support of
$\alpha_2$. (In this discussion, following the notation in Brin's
papers \cite{BrinU} and \cite{BrinEG}, elements of $\ploi$ act on the
right on $I$, and the support of any particular element of $\ploi$ is
the open set of points in $I$ that are moved by that element.)  In
particular, any two distinct conjugates of $\alpha_2$ by powers of
$\alpha_1$ commute with each other, since their supports will be
disjoint in the interval $I$.  Now, if we consider any element
$h\in\langle \alpha_1,\alpha_2\rangle$, it is a standard algebraic
fact that we can write $h$ as a product of the form a power of
$\alpha_1$ followed by a product of conjugates of $\alpha_2$ (and
$\alpha_2^{-1}$) by various powers of $\alpha_1$.  In particular, the
element $\alpha_1$ generates a group $T$ isomorphic to $\Z$ which acts
on the normal subgroup $B$ of $\langle \alpha_1,\alpha_2\rangle$
generated by the conjugates of $\alpha_2$ by different powers of
$\alpha_1$.  Since these conjugates all have disjoint support, the
group $B$ is isomorphic to a direct sum of copies of the integers.
The following chain of isomorphisms should now make sense:
\[
\langle \alpha_1,\alpha_2\rangle\cong (\bigoplus_{i\in\Z}\Z)\rtimes\Z\cong\Z\wr\Z\cong W_2.
\]
Note before we move on that $\alpha_1$ and $\alpha_2$ are both
elements of Thompson's group $F$.

All of the $W_i$ can be realized in Thompson's group $F$ in an
entirely similar way, using one-orbital generators; take $\alpha_i$ to
be the conjugate $\alpha_i = \alpha_{i-1}^s = s^{-1}\alpha_{i-1}s$,
where $s$ is the ``shrinking function'' $s(x) =
\frac{1}{4}x+\frac{1}{4}$, and the conjugation takes place in
$PL_o(\R)$, where we replace $\alpha_1$ by the function in $PL_o(\R)$
which behaves as the identity outside of the unit interval $[0,1]$ for
this inductive definition.  Given $i\in\N$, $i\geq 2$, the support of
$\alpha_i$ is contained in a single fundamental domain of
$\alpha_{i-1}$, so that given $j\in\N$, $j\geq 1$, we have $W_j\cong
\langle \alpha_1,\alpha_2,\ldots,\alpha_j\rangle$.  Since conjugating
any element of $F$ by $s$ will still produce an element of $F$, we see
that the $W_i$'s can all be realized in $F$.

\subsubsection{Orbitals, Towers, and Derived Length}
\label{orbitalsTowersLength}

Much of the language of this section is motivated by thinking of
subgroups of $\ploi$ as permutation groups acting on the set $I =
[0,1]$.

If $H$ is a subgroup of $\ploi$, then its support naturally falls into
a collection of disjoint, open intervals.  Each such we will call an
\emph{orbital}\index{orbital!group} of the group $H$.  Given an
element $\gamma\in H$, the group $\langle\gamma\rangle$ has its own
orbitals, which are the connected components of the support of
$\gamma$.  We will use the symbol $o\gamma$ to denote the number of
orbitals of the group $\langle\gamma\rangle$.  Now, given such an
interval $A = (a,b)\subset [0,1]$, we call $A$ an
\emph{orbital}\index{orbital!element} of $\gamma$, we may also refer
to $A$ as an ``element orbital.''  Since $H$ now must have infinitely
many elements that also have $A$ as an orbital, we typically use
\emph{signed orbitals}, which are pairs of the form $(A,\gamma)$ to
indicate not only which element orbital we are considering, but also
the specific element which is ``owning'' that orbital in our
discussion.  Given a signed orbital $s = (A,\gamma)$, we call $A$ the
\emph{orbital} of $s$ and we call $\gamma$ the \emph{signature} of
$s$.  Given a set $X$ of signed orbitals, we will use $O_X$ to denote
the set of element orbitals that are orbitals of elements of $X$, and
we will use $S_X$ to denote the set of signatures that are signatures
of elements of $X$.

It is a standard fact that $\ploi$ can be totally ordered, and the set
of open intervals in $I$ is partially ordered by inclusion, so the set
of signed orbitals is a poset under the induced lexicographic order.
The ordering on $\ploi$ will not play a role here, except by enabling
simplified language in the upcoming definition.

In the example groups $W_i$, each of the $\alpha_k$ had an orbital
whose closure was fully contained in the orbital of $\alpha_{k-1}$,
whenever $k>1$.  The ease with which we can analyze the groups
generated by the elements $\alpha_k$, which have their orbitals
arranged in such a nice ``stack'', motivates the following.

Given a group $G\leq\ploi$ and a set $T$ of signed orbitals of $G$, we
will say $T$ is a \emph{tower of $G$}\index{tower!of group} if $T$
satisfies the following properties:

\be

\item $T$ is a chain in the partial order on the signed orbitals of $G$.

\item For any orbital $A\in O_T$, $T$ has exactly one element of the
  form $(A,g)$.

\ee  

Let $T$ be a tower of $\ploi$.  We will call the cardinality of $T$
its \emph{height}\index{tower!height}, using the simple descriptive
\emph{infinite}\index{tower!infinite} if $T$ has an infinite
cardinality.

We are now in a position to approach the main geometric result of
\cite{bpgsc}.  We will say that a group has depth $n\in\N$ if and only
if we can find towers of height $n$, but no towers of height $n+1$, in
the group.  The main result of \cite{bpgsc} is as follows:

\bt
\label{geoClassification}

Suppose $G$ is a subgroup of $\ploi$ and $n\in\N$.  G
is solvable with derived length $n$ if and only if $G$ has depth $n$.

\et

This geometric result is the cornerstone upon which the results of
this paper are built.

\section{Classification of solvable subgroups in $\ploi$\label{solveClassificationSection}}
In this section we will pursue the algebraic classification of the
solvable subgroups in $\ploi$.  One direction of the classification is
mostly algebraic, and requires less knowledge of the terminology of
$\ploi$.  We will engage in that direction first.

\subsection{The class $\ws$\label{WS}}
Recall that $\ws$ represents the smallest non-empty class of
isomorphism classes of groups which is closed under the following
three operations (again, given as operations on groups).

\be
\item Restricted wreath product with $\Z$.
\item Bounded direct sum.
\item Taking subgroups.
\ee 

It turns out that we can begin to understand $\ws$ more deeply via a
close study of the class $\mm =
\left\{G_i\,|\, i\in\N\right\}$ defined in the introduction.
Let us gather some facts about the groups $G_k$. (Below, we will prove
the second point in Remark \ref{rkGn}, which corresponds to the fourth
point here.)

\bl 
\label{MFacts}

\be

\item If $F_0$, $F_1$, $H_0$, and $H_1$ are groups, where $F_0\leq
F_1$ and $H_0\leq H_1$, then $F_0\wr H_0\leq F_1\wr H_1$.
\item For any $m$ and $n\in\N$, with $m< n$, $G_m$ embeds
  as a normal subgroup of $G_n$.
\item For any group $G$ with derived length $n$, the groups $G\wr\Z$
and $\bigoplus_{i\in\Z}(G\wr\Z)$ have derived length $n+1$.
\item Given any $n\in\N$, $G_n$ has derived length $n$.

\ee
\el

pf: The first point is immediate by examining the following chain
of subgroup inclusions, where the inclusions are based on the
underlying sets.
\[
\begin{array}{l}
F_0\wr H_0 = \left\{((f_0)_{h\in H_0},h_0)\,|\,(f_0)_{h\in
H_0}\in\bigoplus_{a\in H_0}F_0, h_0\in H_0\right\}
\leq \\
F_1\wr H_0 = \left\{((f_1)_{h\in H_0},h_0)\,|\,(f_1)_{h\in
H_0}\in\bigoplus_{a\in H_0}F_1, h_0\in H_0\right\}
\leq \\
F_1\wr H_1 = \left\{((f_1)_{h\in H_1},h_1)\,|\,(f_1)_{h\in
H_1}\in\bigoplus_{a\in H_1}F_1, h_1\in H_1\right\}
\end{array}
\]
To see the second point, we will demonstrate an embedding of $G_{n-1}$
into $G_n$, and thus inductively define an embedding of $G_m$ into
$G_n$, for any non-negative integers $m <n$.  Note that there are many
copies of $G_{n-1}$ in $G_n$, but we are particularly interested in
the one given in the next paragraph, which is the copy that we use to
inductively define our particular copy of $G_m$ in $G_n$.  The
normality of this embedded copy of $G_m$ in $G_n$ follows easily from
the theory of group actions, which can be checked in section
\ref{realizingGn} below where we realize $G_n$ in $\ploi$.  (Let $X$ be
the set of left hand endpoints of the components of the support of
$G_m$ in $I$.  $X$ is acted upon by $G_n$, and the kernal of this
action is $G_m$.)  A second proof is by noting that the embedded copy
of $G_{n-1}$ in $G_n$ that we demonstrate in the paragraph below is
characteristic in $G_n$, which proof can be carried out with the help
of the geometric tools established in the proof of Theorem
\ref{geoClassification}.  We will not use the normality of our
embedded copy of $G_m$ in $G_n$ later.

Now let us describe our particular embedded copy of $G_{n-1}$ in
$G_n$.  First, identify the base group of $G_{n-1}\wr\Z$ with
$G_{n-1}$ using the fact that $\bigoplus_{i\in\Z}G_{n-1}\cong
G_{n-1}$, now since the direct sum of the base groups of the
$G_{n-1}\wr\Z$ summands in the definition of $G_n$ is also a subgroup
of $G_n$, we see that $\bigoplus_{i\in\Z}G_{n-1}$ is a subgroup of
$G_n$.  But now again, this last direct sum is isomorphic with
$G_{n-1}$, so that $G_{n-1}$ (as embedded here as the direct sum of
the base groups of the $G_{n-1}\wr\Z$ summands of $G_n$) is a subgroup
of $G_n$.

For the third point, it follows from Neumann \cite{NeumannW} (Theorem
4.1 and Corollary 4.5) that if $G=A\wr B$, then $G'$ is contained in a
sum of copies of $A$, and also that $G'$ surjects onto $A$ (these
facts are under the condition that $B$ is abelian).  Both facts are
easy exercises in our situation with $B=\Z$.  This immediately implies
the third point.

The fourth point follows directly from the third.

\qquad$\diamond$

Now let us examine $\ws$ for a short time.  We give (modulo work to
come later) a characterization of $\ws$ that completes one direction
of Theorem \ref{solveClassification}.  The key result (Lemma
\ref{threeClassifications} below) will not be used in the rest of the
paper.

Since $\ws$ is nonempty and closed under taking subgroups, the trivial
group $1$ is a group in $\ws$.  Let $\cg$ represent the class of all
isomorphism classes of groups.  Define $\pc:\pow{\cg}\to\pow\cg$ to be
the closure operator that takes a class $X$ of isomorphism classes of
groups and computes the smallest class of isomorphism classes of
groups which contains $X$ and is closed under the isomorphism class
operations induced from the group operations of restricted wreath
product with $\Z$ and bounded direct sum.  Since each $G_i$ is
obtained by applying a finite sequence of bounded direct sums and
wreath products with $\Z$ to the trivial group $1$, we see that
$\mathscr{M}\subset \left\{1\right\}\pc$.  By the definition of $\ws$,
it is immediate that $\left\{1\right\}\pc\subset \ws$.  In particular,
we have:

\[
\mathscr{M} \subset \left\{1\right\}\pc\subset \ws
\] 

Now let us consider three operators
\[
\begin{array}{l}
\sc:\pow{\cg}\to\pow{\cg},\\
\wc:\pow{\cg}\to\pow{\cg},\\
\bsc:\pow{\cg}\to\pow{\cg},
\end{array}
\]
that represent taking closure under the operations of taking
subgroups, taking restricted wreath products with $\Z$, and building
bounded direct sums, respectively.

If $\Gamma$ is the set of all finite length words of the form
$(\wc\bsc)^k$ or $(\bsc\wc)^k$ where $k\in\N$, and $X\in\pow{\cg}$,
then the union $\Upsilon_{X} = \cup_{\gamma\in\Gamma}X\gamma$ is the
smallest closed set containing $X$ that is closed under both
operations of taking bounded direct sums and wreath products with
$\Z$, so that $\Upsilon_{\left\{1\right\}} = \left\{1\right\}\pc$.

\bl 
\label{BWinM}

Let $m$ be a non-negative integer.  If $G$ is a group in
$\left\{1\right\}\pc$ with derived length $m$, then $G$ is isomorphic
to a subgroup of $G_m$.  
\el

pf: We can prove this by inducting on the derived length of $G$.  If
$G$ has derived length $0$, then $G$ is the trivial group so $G =
G_0$.  Let $n\in\N$, and suppose $G$ has derived length $n$ and that
for any group $H$ in $\left\{1\right\}\pc$ which has derived length
$m$ where $0\leq m<n$, we know that $H\leq G_m$.  Now, there is a
$j\in\N$ so that $G\in\left\{1\right\}(\bsc\wc)^j$, or
$G\in\left\{1\right\}(\wc\bsc)^j$.  Note that if
$G\in\left\{1\right\}(\wc\bsc)^j$ then
$G\in\left\{1\right\}(\bsc\wc)^{j+1}$, so we will assume that
$G\in\left\{1\right\}(\bsc\wc)^k$ for some minimal non-negative
integer $k$.  There are now two cases:

\be

\item $G\in\left\{1\right\}(\bsc\wc)^{k-1}\bsc$ but $G$ is not in
$\left\{1\right\}(\bsc\wc)^{k-1}$.

In this case, The last operation required to build $G$ was a bounded
direct sum of groups, all of which groups have derived lengths
necessarily less than or equal to $n$ (note here that a finite
sequence of bounded direct sums is isomorphic to a bounded direct
sum).  The summands of this bounded direct sum are all groups in
$\left\{1\right\}(\bsc\wc)^{k-1}$.  We will argue that each of these
groups is actually a subgroup of $G_n$, and therefore, since a
countable or finite direct sum of groups isomorphic to $G_n$ is
actually isomophic to $G_n$, we will have finished this case.

Let $H$ be a summand of the final bounded direct sum which created
$G$, and assume that the derived length of $H$ is actually $n$ (at
least one summand must have this derived length), and that
$H\in\left\{1\right\}(\bsc\wc)^{k-1}$.  If $H$ is actually in
$\left\{1\right\}(\bsc\wc)^{k-2}\bsc$, then we can replace $H$
inductively by a summand (with derived length $n$) of the last bounded
sum operation used to create $H$, so that there is a $t\in\N$ so that
$H$ is now an element of $\left\{1\right\}(\bsc\wc)^t$, but not an
element of $\left\{1\right\}(\bsc\wc)^{t-1}\bsc$.  In particular, $H$
is the result of applying $s$ wreath products with $\Z$ to a group
$H^*$ in $\left\{1\right\}(\bsc\wc)^{t-1}\bsc$ for some $s\in\N$.
$H^*$ is therefore a group in $\left\{1\right\}\pc$ with derived
length $n-s$, and therefore $H^*$ is a subgroup of $G_{n-s}$.  But
note that for any integer $p$ we have that $G_p\wr \Z\leq G_{p+1}$, so
abusing notation (the parenthesies below should accumulate on the left),
we see that $H = H^*(\wr\Z)^s\leq G_{n-s}(\wr\Z)^s\leq
G_{n-s+1}(\wr\Z)^{s-1}\leq \ldots \leq G_n$.

If $H$ is a summand of $G$ with derived length $m$ where $m<n$, Then by
our induction hypothesis, $H\leq G_m$, but $G_m\leq G_n$, so $H\leq
G_n$.

\item $G\in\left\{1\right\}(\bsc\wc)^k$ but $G$ is not in
$\left\{1\right\}(\bsc\wc)^{k-1}\bsc$.

This case is entirely similar to the last case, except that we already
know that $G$ is the result of applying $s$ wreath products with $\Z$
to a group $H$ in $\left\{1\right\}(\bsc\wc)^{k-1}\bsc$ for some
positive integer $s$.  The derived length of $H$ must be $n-s$, and
therefore $H\leq G_{n-s}$ so that $G\leq G_n$ as in the penultimate
paragraph of the previous case.

\ee 
\qquad$\diamond$

We need one last technical lemma before we can complete our
exploration of the class $\ws$:

\bl 
\label{shallowGroupEmbeddings}

Suppose $m$ and  $n\in\N$, with $m\leq n$. If $H\leq G_n$ and $H$ has
derived length $m$, then $H$ is isomorphic to a subgroup of $G_m$.

\el

pf: This follows from Corollary \ref{solveEmbeddings} below, and the
fact that we can realize the groups $\mm$ in $\ploi$ (see the next
subsection).   \qquad$\diamond$

Our arguments after this section do not rely on the
classification of $\ws$ given in the next lemma.

Finally, we have a nice description of $\ws$.  Note that the following
lemma implies Theorem \ref{RandM}.
\bl
\label{threeClassifications}

$\ws = \left\{1\right\}\pc\sc = \mathscr{M}\sc$
\el

Pf: 
We have already shown that $\left\{1\right\}\pc \subset
\mathscr{M}\sc$, so we know that $\left\{1\right\}\pc\sc \subset
\mathscr{M}\sc$, and by definition, $\mathscr{M}\subset
\left\{1\right\}\pc$, so that $\mathscr{M}\sc\subset
\left\{1\right\}\pc\sc$.  In particular, $\mathscr{M}\sc =
\left\{1\right\}\pc\sc$.

We will now show that $\mathscr{M}\sc\bsc = \mathscr{M} \sc$ (implying
that $\mathscr{M}\sc$ is already closed under the operation of taking
bounded direct sums) and that $\mathscr{M} \sc\wc = \mathscr{M} \sc$
(implying that $\mathscr{M}\sc$ is already closed under taking wreath
products with $\Z$), so that we can conclude that $\mathscr{M} \sc =
\ws$.

To see that $\mathscr{M}\sc\bsc = \mathscr{M}\sc$, let $G\in
\mathscr{M}\sc\bsc$.  There is an $M\in\N$ so that we can write
$G=\bigoplus_{i\in\Z}H_i$, where each $H_i$ has derived length bounded
above by $M$.  Now, each $H_i\leq G_M$ (by the Lemma
\ref{shallowGroupEmbeddings}), so we see that $G \leq
\bigoplus_{i\in\Z}G_M$, hence $G\in \mathscr{M}\sc$.

To see that $\mathscr{M}\sc\wc = \mathscr{M}\sc$, we note that if
$G\in\mathscr{M}\sc\wc$, then either $G\in\mathscr{M}\sc$ or we can
write $G = ((\ldots((H\wr \Z)\wr\Z)\ldots)\wr\Z$, where
$H\in\mathscr{M}\sc$, and where there are $k$ wreath products with
$\Z$, for some $k\in\N$.  In the first case we are done.  In the
second case $H \leq G_M$ for some non-negative integer $M$.  But now
$G=((\ldots(H\wr\Z)\wr\Z)\wr\ldots)\wr\Z\leq
((\ldots(G_M\wr\Z)\wr\Z)\wr\ldots)\wr\Z=J\in \left\{1\right\}\pc$,
where $J$ has derived length $M+k$, so that $G\leq G_{M+k}$.  Finally
we have that $G\in \mm\sc$.  
\qquad$\diamond$

\subsection{Realizing $\ws$ in $\ploi$ and $F$\label{realizingGn}}
Here we will explain how we can realize the groups in $\mathscr{M}$
inside Thompson's group $F$, as realized in $\ploi$.  This will prove
one half of Theorem \ref{solveClassification}.  To realize the other
direction of Theorem \ref{solveClassification}, we will need to take a
detour through the geometric definitions of $\ploi$.

The elements $\alpha_1$, $\alpha_2\in\ploi$ defined in the
introduction will play a major role here, as will the shrinking
conjugator $s$.

First, observe that we can realize $G_1$ fairly simply.  Let
$\beta_0$ = $\alpha_2$, and define $\beta_k$ for each $k\in\Z$ as
$\beta_0^{\alpha_1^k}$.  $G_1$ is immediately isomorphic to $\langle
\beta_k|k\in\Z\rangle$, and $G_1$ has been realized in Thompson's group $F$.

Now we will show that given any group $H$ which has been realized as a
subgroup of $\ploi$, we can realize $\bigoplus_{i\in\Z}(H\wr \Z)$ as a
subgroup of $\ploi$.  Firstly, conjugate the elements of $H$ by $s$
twice, to create a new group $H_0$ isomorphic with $H$.  The supports
of all of the elements of $H_0$ are contained in the set
$(\frac{5}{16},\frac{3}{8})$, which is contained in a single
fundamental domain of $\alpha_2 = \beta_0$.  At this juncture, the
group generated by $H_0$ and $\beta_0$ is isomophic to $H\wr\Z$, where
$\beta_0$ is the generator of the top $\Z$ factor.  In particular, we
can realize $H\wr\Z$ in $\ploi$.  But now $\bigoplus_{i\in\Z}(H\wr\Z)$
is the base group of $(H\wr\Z)\wr\Z$, which we can also realize by
repeating the previous procedure, so we are done.

  Observe that the shrinking map conjugates elements of Thompson's
group $F$ into Thompson's group $F$, so that if $H$ is realized as a
subgroup of $F$, then this construction of
$\bigoplus_{i\in\Z}(H\wr\Z)$ also produces a subgroup of $F$.  In
particular, we see inductively that each group $G_n$ can now be
realized in Thompson's group $F$, and therefore that each group of $\ws$ can be
realized in Thompson's group $F$.

\section{{\boldmath $\ploi$}}
We will now build required terminology for working with subgroups
of $\ploi$.  We will use notation similar to that in Brin's paper
$\cite{BrinU}$ on the ubiquitous nature of Thompson's group $F$ in
subgroups of $\ploi$.  Some portions of the text here were lifted
from \cite{bpgsc}.

We note that the set of points moved by an element $h$ of $\ploi$ is
open, by the continuity of elements of $\ploi$.  But then the support
of $H$, for any subgroup $H\leq\ploi$, is a countable union of
pairwise disjoint open intervals in $(0,1)$.  Let the collection
$\mathscr{O}_H$ always denote the countable set of open, pairwise
disjoint intervals of the support of $H$.  Recall from section
\ref{orbitalsTowersLength} that we call these intervals the
\emph{orbitals of $H$}\index{orbital!group}. There is a natural total
order on $\mathscr{O}_H$, where if $A$, $B\in\mathscr{O}_H$, where
$A\neq B$, we will say $A<B$ or \emph{$A$ is to the left of
$B$}\index{orbitals!natural ordering of} if and only if given any
$x\in{}A$ and $y\in{}B$, we have $x<y$ under the natural order induced
by $I\subset\R$.  Since $A$ and $B$ are disjoint, connected subsets of
$\R$, this definition does not depend on the choices of $x$ and $y$.

If the collection $\mathscr{O}_H$ is finite, we may speak of the
``first'' orbital, or ``second'' orbital, etc., where the first
orbital is the leftmost orbital under the definition given above, the
second orbital is the orbital to the left of all other orbitals in
$\mathscr{O}_H$ other than itself and the first orbital, and so on.

Given an open interval $A =(a,b)\subset\R$, where $a<b$, we will refer
to $a$ as the \emph{leading end of $A$}\index{end!leading}, and to $b$
as the \emph{trailing end of $A$}\index{end!trailing}.  If the
interval is an orbital of some group $H\in\ploi$, we will refer to
the \emph{ends of the orbital}\index{orbital!end} in the same fashion.

If $h\in{}H$ and $x\in{}\supp(h)$, we will say that \emph{$h$ moves $x$
to the left}\index{movement!left} if $xh<x$, and we will say that
\emph{$h$ moves $x$ to the right}\index{movement!right} if $xh>x$.
Furthermore, we will say that $x\in{}I$ is a \emph{breakpoint for
$h$}\index{breakpoint} if the left and right derivatives of $h$ exist
at $x$, but are not equal.  We recall that by definition, $h$ will
admit only finitely many breakpoints.  If $\mathscr{B}_h$ represents
the set of breakpoints of the element $h$, then
$(0,1)\backslash\mathscr{B}_h$ is a finite collection of open
intervals, which we will call \emph{affine components of
$h$}\index{component!affine}, which admit a natural ``left to right''
ordering as before.  We shall therefore refer to the ``first''
affine component of $h$, or the ``second'' affine component of $h$,
etc.  We sometimes will refer to the first affine component of $h$ as
the \emph{leading affine component of $h$}\index{component!leading
affine}, and to the last affine component of the domain of $h$ as the
\emph{trailing affine component of $h$}\index{component!trailing
affine}.

The following are some useful remarks (whose proofs are left to the
reader), that are mostly standard in the theory of $\ploi$.

\brk
\label{finiteOrbitals}
\label{transitiveElementOrbital}
\be

\item If $A$ is an orbital for $h\in H$, then either $xh>x$ for all
points $x$ in $A$, or $xh<x$ for all points $x$ in $A$.

\item Any element $h\in\ploi$ has only finitely many
orbitals.

\item If $h\in\ploi$ and $A=(a,b)$ is an orbital of $h$, then given any
$\epsilon>0$ and $x$ in $A$, there is an integer $n$ so that
$|xh^{-n}-a|<\epsilon$ and $|xh^n-b|<\epsilon$.
\ee
\erk

Given an orbital $A$ of $H$ we say that \emph{$h$ realizes an end of
$A$} if some orbital of $h$ lies entirely in $A$ and shares an
end with $A$.  Note that Brin uses the word ``Approaches'' for this
concept in \cite{BrinU}, but we will use ``Approaches'' to also
indicate the direction in which $h$ moves points, as follows: we will
say that \emph{$h$ approaches the end $a$ of $A$ in $A$} if $h$ has an
orbital $B$ where $B\subset A$ and $B$ has end $a$, and $h$
moves points in $B$ towards $a$.  In particular, $h$ realizes $a$ in $A$
and $h$ moves points in its relevant orbital towards $a$.  If $A$ is
an orbital for $H$ then we say that \emph{$h\in H$ realizes $A$} if
$A$ is also an orbital for $h$.

If $g$ and $h$ are elements of $\ploi$ and there is an interval $A =
(a,b)\subset I$ so that both $g$ and $h$ have $A$ as an orbital, then
we will say that $g$ and $h$ \emph{share the orbital $A$}.

\subsection{Conjugation and Transition Chains}
All statements in this section will be either fairly standard and
given without proof, or are proven (or are immediate from proofs) in
\cite{bpgsc}.  First, we will give some of the standard facts.

Let $g$, $h\in\ploi$ and let $k = g^h = h^{-1}gh$.  Suppose that
$\mathscr{O}_g=\bSeq{A}{i}{1}{n}$ are the $n$ orbitals of $g$ in left
to right order, where $n\in\N$ and $i\in\left\{1,2,\ldots,n\right\}$.
Define
\[
A_i^*=\left\{x\in{}I|xh^{-1}\in{}A_i\right\}=A_ih
\]
for all $i\in\left\{1,2,\ldots,n\right\}$.

\bl 

$o_k = o_g = n$ and collection $\bSeq{A^*}{i}{1}{n}$ is the ordered
set of orbitals of $g^h$ in left to right order.  

\el

In the setting of the above lemma, we will say that the $A_ih$  are the
\emph{induced orbitals}\index{orbital!induced} of $k$ from $g$ by the
action of $h$.  We might also say that the orbitals of $k$ are induced
from the orbitals of $g$ by the action of $h$.

The following is worth pointing out:
\brk 
\label{breakpoints}
Suppose $g$, $h\in\ploi$ and $f = gh$.  If $b$ is a breakpoint of $f$
then $b$ is a breakpoint of $g$ or $bg$ is a breakpoint of $h$.
\erk

Finally, our last well known fact.

\bl
\label{transitiveOrbital}
If $H\leq\ploi$ and $A=(a,b)$ is an orbital for $H$, then given any points
$c$, $d\in A$, with $c<d$, there is an element $g\in{}H$ so that
\mbox{$cg>d$}.
\el

For the rest of the section, we give definitions and basic
statements introduced in \cite{bpgsc}.

The following definiton is a special case of a more general
definition.  In this paper, we will only need objects of the more
restricted type. Let $A=(a,b)$ and $B=(c,d)$ be two intervals
contained in $I$, and let $\alpha$, $\beta\in G$ for some subgroup
$G\leq \ploi$ so that the set
$\mathscr{C}=\left\{(A,\alpha),\,(B,\beta)\right\}$ is a set of signed
orbitals for $G$.  We will say $\mathscr{C}$ is a \emph{transition
chain of length two} for $G$ if $a<c<b<d$.  The existence of a
transition chains of length two in a group allows the
creation of complex dynamics in the group action on the unit interval,
as one can repeatedly act on the interval with one element to move a
point near to an edge of an element orbital, and then use the ``next''
element to move the point out of its initial element orbital.

\subsection {Notes on Towers and Orbitals}
In order to finish the proofs of the main theorems of this paper, we
will need a deeper understanding of towers.  Recall that a tower $T$ is a
set of signed orbitals satisfying the following two conditions:
\be
\item $T$ is a chain in the poset of signed orbitals associated with
  $\ploi$, and 
\item for any orbital $A\in O_T$, $T$ has exactly one
  element of the form $(A,g)$.  
\ee 
We will refine this definition, following the related section 2.3 in
  \cite{bpgsc}.

The notions of depth of groups and height of towers can be
extended to other objects.

Given a group $G\leq \ploi$, and an open interval $A= \subset I$, we
will define the \emph{depth of $A$ in $G$}\index{orbital!depth} to be
the supremum of the heights of the finite towers which have their
smallest element having the form $(A,h)$ for some element $h\in G$.
In particular, if $A$ is not an orbital for an element of $G$, then
the depth of $A$ will be zero.  Symmetrically, we define the
\emph{height of $A$ in $G$}\index{orbital!height} to be the supremum
of the heights of the finite towers which have their largest element
having the form $(A,h)$ for some element $h\in G$.

Towers, as defined, are easy to find, but can be difficult to work
with.  For an arbitrary tower $T$, there are no guarantees about how
other orbitals of signatures of the elements of $T$ cooperate with the
orbitals of the tower.  We say a tower $T$ is an \emph{exemplary
tower}\index{tower!exemplary} if the following two additional
properties hold: \be

\item Whenever $(A,g)$, $(B,h)\in T$ then $(A,g)\leq(B,h)$ implies the
orbitals of $g$ are disjoint from the set of ends of the orbital $B$.

\item Whenever $(A,g)$, $(B,h)\in T$ then
$(A,g)\leq(B,h)$ implies no orbital of $g$ in $B$ shares an end with $B$.

\ee

\subsection{Thompson's group $F$ and balanced subgroups of $\ploi$}
Brin showed in \cite{BrinU} the following theorem:

\bt [Ubiquitous F]
\label{UbiquitousF}  

If a group $H \leq \ploi{}$ has an orbital $A$ so that some element
$h\in{}H$ realizes one end of $A$, but not the other, then $H$ will
contain a subgroup isomorphic to Thompson's group $F$.
\et

  We will say that an orbital $A$ of a group $H \leq \ploi{}$ is
\emph{imbalanced}\index{orbital!imbalanced} if some element $h \in H$
realizes one end of $A$, but not the other, and we will say $A$ is
\emph{balanced}\index{orbital!balanced} if whenever an element $h \in
H$ realizes one end of $A$, then $h$ also realizes the other end of
$A$ (note that $h$ might do this with two distinct orbitals).
Extrapolating, given a group $H \leq \ploi{}$, we will say that
\emph{$H$ is balanced}\index{group!balanced} if given any subgroup
$G\leq H$, and any orbital $A$ of $G$, every element of $G$ which
realizes one end of $A$ also realizes the other end of $A$.
Informally, $H$ has no subgroup $G$ which has an orbital that is
``heavy'' on one side.  In the case where $H$ has a subgroup $G$ with
an imbalanced orbital, then we will say that \emph{$H$ is
imbalanced}\index{group!imbalanced}.

\brk 

If $H\leq\ploi$ and $H$ is imbalanced, then $H$ has a subgroup
isomorphic to Thompson's group $F$.

\erk

Since $F'$ is a non-abelian simple group (\cite{CFP},
Theorem 4.3), $F$ is not solvable.  Thus imbalanced groups are not
solvable.

The dynamics of balanced groups are much easier to understand than
those of imbalanced groups.  One indicator of this is that imbalanced
groups admit transition chains of length two, and groups which admit
transition chains of length two admit infinite towers.

The following lemma sums up what we have learned, from a utilitarian
perspective.  It is an easy consequence of the results in \cite{bpgsc}.

\bl 
\label{solvableQualities}

If $G$ is a solvable subgroup of $\ploi$ then

\be

\item $G$ is balanced,
\item $G$ does not admit transition chains of length two, and
\item all towers for $G$ are exemplary.  

\ee
\el

\subsubsection{A useful homomorphism}
\label{phi}
The homomorphism in this section was known to Brin and Squier during
the research that led to the paper \cite{picric}.

Let us suppose that $H \leq \ploi$ and $A=(a,b)$ is an orbital of $H$.
To simplify the arguments for now, let us suppose that $A$ is the only
orbital of $H$.  We can define a map $\phi:H \to\R\times\R$ defined by
$h\mapsto(h_a,h_b)$ where $h_a=\ln(h_+'(x))$ and $h_b=\ln(h_-'(x))$.
Ie., we take the logs of the slopes of $h$ at the ends of $A$.  Since
$h$ is a piecewise-linear orientation-preserving homeomorphism of $I$,
we see that the derivatives exist and are positive, and so $\phi$ is
well defined for all $h\in H$.  If $h$ does not realize $a$
(resp. $b$) then we see that $h$ behaves as the identity near $a$
($b$) in $A$, and so $h_a = \ln(1)=0$ ($h_b = 0$).  If $h,g \in H$
then $hg\phi=h\phi+g\phi$ in $\R\times\R$ by the chain rule. In
particular, we see the following remark:

\brk
$\phi$ is a homomorphism of groups.
\erk

Now the image of $\phi$ is quite interesting, it carries a small
amount of the complexity of $H$, but still enough to allow us to find
out if $A$ is an imbalanced orbital.  The following straightforward
lemma is left to the reader.

\bl
The orbital $A$ is imbalanced if and only if $Im(\phi)$ contains an
element of the form $(\alpha,0)$ or $(0,\alpha)$ where $\alpha \neq
0$.
\el 

This next technical lemma will help with the lemma that follows it:
\btl 
\label{longLeadSlope}

Suppose $H\leq\ploi$, $H$ has an orbital $A=(a,b)$, and that $H$
has a sequence of elements $(g_n)_{n = 0}^{\infty}$ in $H$ which
satisfies the properties below.

\be

\item For each $i\in\N$, the lead slope of $g_{i+1}$ in $A$ is less than the
lead slope of $g_i$ in $A$.

\item Given any real number $q>1$, there is an $i\in\N$ so that the
lead slope of $g_i$ in $A$ is $p$ where $1<p<q$.

\ee

Then there is $c\in(a,b)$ so that given any real number $s>1$, $H$ has an
element $\alpha$ which has an affine component $\Gamma$ containing
$(a,c)$, and with slope $r$ on $\Gamma$ where $1<r<s$.

\etl

pf:

To simplify the language of this argument, we will restrict our
attention to the orbital $A$, treating it as the domain of the
elements of $H$, so that the phrase ``The first affine component of
[$h\in H$]'' will really mean the open interval $(a,u)$ where $h$ has
an affine component of the form $(v,u)$ where $v \leq a<u$.  We will also
refer to this as the ``leading (or lead) affine component of $h$''.
Similarly, we will refer to the slope of $h$ on its leading affine
component as the ``lead slope of $h$.''

Note that the second condition on $(g_n)_{n=0}^{\infty}$ implies that
every element of $(g_n)_{n=0}^{\infty}$ has lead slope greater than
one.

For each $i\in\N$, let $(a,b_i)$ be the first affine component of
$g_i$, and let $s_i$ represent the lead slope of $g_i$.

We are now in a position to define a new sequence of
functions $\pSeq{h}{i}{0}{\infty}$ which satisfy the following conditions:

\be

\item For each $i\in\N$, the lead slope of $h_i$ is $s_i$.

\item For each $i\in\N$, the leading affine component of $h_i$
contains $(a,b_0)$.

\ee

Given such a sequence $(h_i)_{i=0}^{\infty}$ and $s>1$, the
hypothesies on the $g_i$ show that there will be an $N\in\N$ so that
for all $n>N$ we will have that $s>s_n>1$.  If we further set $c =
b_0$, then for all $n\in\N$, $h_n$ will have leading affine component
containing $(a,c)$, so we will have finished the proof of the lemma.

For each $i\in\N$, define $n_i$ to be the smallest non-negative
integer so that $b_0g_0^{-n_i}<b_i$, and define $h_i =
g_0^{-n_i}g_ig_0^{n_i}$.  Since $g_0$ moves points to the right on its first
affine component, $n_i$ is well defined, and therefore $h_i$ is well
defined.  We now check that $h_i$ satisfies the two conditions, for
each $i$.  Firstly, we observe that the lead slope of each $h_i$ is
the product of the lead slopes of the elements of the product
$g_0^{-n_i}g_ig_0^{n_i}$, which is $(\frac{1}{s_0})^{n_i}s_is_0^{n_i}
= s_i$.  Secondly, by Remark \ref{breakpoints}, we know that if
$x\in(a,b)$ is a breakpoint of the product $g_0^{-n_1}g_ig_0^{n_i}$,
then the image of $x$ under application of some initial partial
product (possibly empty) must be a breakpoint of the next term of the
overall product.  However, the first breakpoint of $g_0^{-1}$ is $d_0
= b_0g_0>b_0$, and $g_0^{-1}$ moves points left in its first affine
component so if $x\in (a,b_0)$, then $xg_0^{-k}<d_0$ for all
non-negative integers $k$.  In particular, if $x\in(a,b_0)$, then $x$
is in the first affine component of $g_0^{-n_i}$ for any $i\in\N$.
Given an $i\in \N$, the first breakpoint of $g_i$ in $A$ is $b_i$, but
$n_i$ was chosen so that $b_0g_0^{-n_i}<b_i$, so if $x\in(a,b_0)$,
then the image of $x$ under $g_0^{-n_i}$ is in the first affine
component of $g_i$, in particular, $(a,b_0)$ is contained in the first
affine component of $g^{-n_1}g_i$.  Finally, the first breakpoint of
$g_0^{n_1}$ is the image of $b_0$ under $g_0^{-n_i + 1}$, but
$b_0g_0^{-n_i}g_0\geq b_0g_0^{-n_i}g_i$ because the leading slope of $g_i$
is less or equal to the leading slope of $g_0$.  In particular, the whole
interval $(a,b_0)$ is in the first affine component of $h_i$.
\qquad$\diamond$

The following lemma expresses the main point of the section, as it will
enable us to find the remarkable ``controllers''; elements that
control the global behavior of a balanced group on its
orbitals.

\bl
\label{phiImageLemma}

Suppose $H$ is a balanced group with single orbital $A = (a,b)$, and
$\phi$ is the log-slope homomorphism defined before, then
$\phi(H)\cong \mathbb{Z}$ or $\phi(H)$ is trivial in $\rtr$.

\el 

pf: Let $H$ be a balanced subgroup of $\ploi$.  Each element of $H$
either realizes both ends of $A$, or neither.  In particular, the
group homomorphism $\rho_1:\R\times\R\to\R$ which is projection on the
first factor has the property that $ker(\rho_1)\cap{}Im(\phi) =
\{(0,0)\}$, the trivial subgroup of $\R\times\R$.  This implies that
the image of $\phi$ in $\R\times\R$ is isomorphic to the group
$H_1\leq\R$ obtained by considering only the first factors of elements
of $Im(\phi)$.  If no element realizes the ends of $A$, then $H_1 =
\left\{0\right\}$, the trivial (additive) group, and we are done.
Therefore, let us suppose instead that some elements in $H$ realize
the left end of $A$ (and therefore also the right) so that $H_1$
cannot be the trivial subgroup of $\R$.

If $H_1$ is discrete in $\R$ then $H_1$ is either trivial, or
isomorphic to $\Z$, but by assumption, $H_1$ is not the trivial group,
hence in this case $H_1\cong\Z$.  Hence, we shall suppose that $H_1$
is not discrete in $\R$.  In this case, by taking the difference of
two elements in the image which are very near each other, we see that
we can find an element of $H_1$ which is as close to zero as we like.
This implies there are elements of $H$ whose leading slopes are as
close to one as we like, without actually being one.  If $h$ is an
element with leading slope $s\neq 1$, then one of $h$ or $h^{-1}$ has
slope greater than one, since the leading slope of $h$ is $s$, but the
leading slope of $h^{-1}$ is $1/s$ (note that $s$ cannot be zero,
since no element of $\ploi$ has an affine component with slope zero).

Now suppose that $H$ is abelian.  By a result of Brin and Squier in
\cite{picric}, if two elements in $\ploi$ commute and share a common
orbital $W$, then there is an element $w\in\ploi$ which behaves as a
the identity off of $W$ and is a common root of the initial two
elements over $W$.  If two elements have non-disjoint support and
commute, it is easy to see that the intersections of their supports
actually is a set of commonly shared orbitals. Now, let $h$ be some
element of $H$ with leading slope $s>1$.  For each positive integer
$n$, let $g_n$ be an element of $H$ with leading slope $s_n'$ where
$1<s_n'<\frac{n+1}{n}$.  Now for each $g_n$, the pair $h$ and $g_n$
has a common root $h_n$ (on their leading orbital), but infinitely
many of the roots $h_n$ have pairwise distinct leading slopes, since
these slopes are always less than or equal to the slopes of the $g_n$,
and in particular, $h$ must then have infinitely many distinct roots
in $H$ on its leading orbital.  By another result in \cite{picric}, no
element of $\ploi$ has infinitely many distinct roots, so we must
conclude that $H$ is not abelian.

Note that by the details of the previous paragraph, it is easy to construct a
countably infinite sequence of elements $(g_n)_{n = 1}^{\infty}$ in $H$ which
satisfies the properties below: 
\be

\item For each $i\in\N$, the lead slope of $g_{i+1}$ is less than the
lead slope of $g_i$.

\item Given any real number $r>1$, there is an $i\in\N$ so that the
lead slope of $g_i$ is $s$ where $1<s<r$.

\item given $i$, $j\in\N$, we will have $[g_i,g_j]=1$ implies $i = j$.

\ee

Therefore, by Lemma \ref{longLeadSlope} there is $c\in(a,b)$ and
elements of $H$ with lead slopes that are greater than but arbitrarily
close to one, and whose leading affine components contain $(a,c)$.

Let $f$ and $g$ be two elements of $H$ with $h = [f,g]\neq 1$.  The
fixed set of $h$ in $I$ is disconnected, and contains two components of the
form $[0,u']$ and $[v',1]$ for some numbers $u'$, $v'\in (a,b)$.  In
particular, $\inf(\supp(h)) = u'$ and $\sup(\supp(h)) = v'$.

By Lemma \ref{transitiveOrbital} there is $\alpha\in H$ so that
$v'\alpha<c$, so that $j=h^{\alpha}$ has all of its orbitals inside
$(a,c)$.  Let $u = \inf(\supp(j))$ and $v = \sup(\supp(j))$.  In
particular, $v<c$.

  We will now perturb $j$ slightly via a conjugation which will move
the orbitals of $j$ to the right by a distance less than $L$, so that
$j$ and the new element together will generate a group with an
imbalanced orbital.  

Suppose $L>0$ is smaller than two particular lengths.  The first
length is the length of the second component of the fixed set of $j$
in $A$ which has non-zero length (note, this component might just be
$[v,b)$, if $j$ has only two such), and the second length is the
length of the first orbital of $j$.

Choose an element $\beta\in H$ whose leading affine component in $A$
contains $(a,c)$ and whose lead slope is greater than one, but so near
one that no point in $(a,c)$ will move to the right a distance greater
than $L$.  Now the elements $j$ and $j^{\beta}$ will generate a group
$G$ with leading orbital $(u,w)$ where $j$ realizes $u$ and possibly
$w$ (if the right ends of the appropriate orbitals of $j$ and
$j^{\beta}$ are aligned), while $j^{\beta}$ will achieve $w$ but not
$u$.  To see this, note that the left end of the first orbital of
$j^{\beta}$ is in the first orbital of $j$, so that $u$ is the left
end of the leading orbital of $G$, and only $j$ realizes it.
Meanwhile, the right end of the first orbital of $j^{\beta}$ is to the
right of the right end of the first orbital of $j$, so that the
leading orbital of $j^{\beta}$ contains the right end of the leading
orbital of $j$.  As we progress to the right, if the solitary fixed
point sets of $j$ and $j^{\beta}$ align before we reach the second
component of the fixed set of $j$ in $A$ with non-zero length, then
the right end of the first orbital of $G$ will be achieved by both $j$
and $j^{\beta}$.  Otherwise, the first orbital of $G$ will extend
rightward into the interior of the second fixed set of $j$ with a
non-empty interior, so that the right end of this orbital of $G$ will
be realized only by $j^{\beta}$.  In all cases, $G$ will be
imbalanced, and hence $H$ will also be imbalanced.  This contradicts
the hypothesies of the lemma, and so we see that $H_1$ must be
discrete in $\R$, and therefore the image of $\phi$ is isomorphic to
$\Z$ or the trivial group in $\R\times\R$.  \qquad$\diamond$

The kernal of the homomorphism $\phi$ is naturally very important as
well, it is the subgroup of $H$ which consists of elements which are
the identity near the ends of $A$.  Typically, we will refer to this
normal subgroup as ${}H^{\!\!\!\!^{^\circ}}$.

\subsubsection{Controllers}
\label{controllers}

A consequence of section \ref{phi} is that the structure of a balanced
group with one orbital is very special.  In this section we will
explore this idea further.

\bl [Balanced Generator Existence]
\label{balancedGenerator}

Suppose that $H$ is a balanced subgroup of $\ploi$ with single orbital
$A$, that there is some element in $H$ which realizes an end of $A$,
and that ${}H^{\!\!\!\!^{^\circ}}$ is the subgroup of $H$ which
consists of all elements in $H$ which are the identity near the ends
of $A$. Then there is an element $g$ of $H$ so that
$H=\left<\left<g\right>,{}H^{\!\!\!\!^{^\circ}}\right>$, where $g$
realizes both ends of $A$.  
\el

pf: 

Let $\Gamma_A$ be the set of elements of $H$ which realize both ends
of $A$.  By our assumptions, $\Gamma_A$ is not empty.  Now observe
that $H = 
\left<\left<\Gamma_A\right>,\,{}H^{\!\!\!\!^{^\circ}}\right>$.

By lemma \ref{phiImageLemma} the image $\phi(H)$ is infinite cyclic in
$\R\times\R$.  Let $\gamma$ be a generator of the image of $\phi$.
Let $g$ be an element of $H$ so that $g\phi =\gamma$.  We observe that
since $\gamma$ is non-trivial in both components, $g$ realizes both
ends of $A$.  Since $\gamma$ generates the image of $\phi$, if
$\hat{g}\in{}\Gamma_A$, then $\hat{g}\phi = g^k\phi$ for some $k\in\Z$.
Hence, $g^{-k}\cdot\hat{g}\in{}H^{\!\!\!\!^{^\circ}}$.  This now implies that
$\Gamma_A\subset \left<g,{}H^{\!\!\!\!^{^\circ}}\right>$, so that
$\left<g,{}H^{\!\!\!\!^{^\circ}}\right> = \left<\Gamma_A,{}H^{\!\!\!\!^{^\circ}}\right> =
H$.\qquad$\diamond$

We will call any element $c$ of a balanced group $H$ with one orbital
$A$, which satisfies the rule $H =
\left<\left<c\right>,\,{}H^{\!\!\!\!^{^\circ}}\right>$, a
\emph{controller}\index{controller} of $H$.  A controller of $H$ is
clearly a special element.

Given a controller $c$ of a balanced group $H$ with one orbital $A$,
we can write any element $h$ of $H$ uniquely in the form
$c^k\cdot g^{\!\!\!^{\circ}}$, where $k$ is some integer, and $g^{\!\!\!^{\circ}}\in{}H^{\!\!\!\!^{^\circ}}$.
We will call this the \emph{$c$-form of $h$}\index{$c$-form}.

We will say that a controller $c$ of the group $H$ is
\emph{consistent}\index{controller!consistent} if and only if its
image $(\alpha,\,\beta) = c\phi$ satisfies the property that
$sign(\alpha) = -sign(\beta)$.  Otherwise we will say the controller
is \emph{inconsistent}\index{controller!inconsistent}.  The idea
behind this definition is that a one orbital controller should be
consistent, since it is either moving points to the right everywhere
on its support, or it is moving points to the left everywhere on its
support.  An inconsistent controller must have a fixed point set, and
has at least one orbital where the controller moves points to the right,
and one orbital where the controller moves points to the left.  It turns
out that a consistent controller actually is a one-orbital element of
$H$.

\bl
\label{cFullSupport}

Suppose $H$ is a balanced subgroup of $\ploi$ and $H$ has unique
orbital $A$.  Further suppose that $H$ has a consistent controller
$c$, then $c$ realizes $A$.

\el

pf:

Since $c$ and $c^{-1}$ are both controllers, and either both satisfy
or both fail the conclusion of the statement of the lemma, we will
assume $c$ moves points to the right on its orbitals near the ends of
$A$.  Suppose $c$ has a non-trivial fixed set $K$ in $A$.  $K$ is
closed and bounded and hence compact.  Let $u = \inf{K}$ and $v =
\sup{K}$, so that $(a,u)$ is the first orbital of $c$ and $(v,b)$ is
the last orbital of $c$.  Now there are points $x\in (a,u)$, and
$y'\in(v,b)$ so that we have $a<x<u\leq v<y'<b$.  By Lemma
\ref{transitiveOrbital} there is an element $g\in H$ so that $xg>y'$.
Writing $g$ in its $c$-form, we have that $g = c^kg^{\!\!\!^{\circ}}$
for some integer $k$ and element
$g^{\!\!\!^{\circ}}\in{}H^{\!\!\!\!^{^\circ}}$.  In particular, the
element $h = gc^{-k}$ is trivial near the ends of $A$, but satisfies
$xh=y>v$. The element $h$ therefore has an orbital $D=(r,s)$ which
spans the fixed set $K$ of $c$.  Suppose $e=\inf(\supp(h))$, so that
$a<e\leq r$.  By Lemma \ref{transitiveElementOrbital} there is a
positive integer $N_1$ so that for any integer $n_1>N_1$ we will have
$r<ec^{n_1}<u$. Suppose $f=\sup(\supp(h))$, so that $s\leq f<b$.  By
Lemma \ref{transitiveElementOrbital} there is a positive integer $N_2$
so that for any integer $n_2>N_2$ we will have $f<sc^{n_2}<b$.  Let $n
= \max(N_1,N_2)$, then, the element $j = h^{(c^n)}$ has its first
orbital starting at some interior point of $(r,s)$, and its orbital
induced from $(r,s)$ has right end $t$ which is strictly to the right
of $f$.  In particular the group $H_1=\left<j,h\right>$ has an orbital
$B = (r,t)$, where we note that $t>s$ by construction.  Now, $h$
realizes the left end of $B$ in $B$, but not the right, hence $H_1$,
and therefore $H$, is imbalanced.  But this contradicts our
assumptions, therefore $c$ must have orbital $A$.\qquad$\diamond$

The following corollary is left to the reader.
\bc 
\label{consistentRealization}
\label{iAntiRealization}
Suppose $H$ is a balanced subgroup of $\ploi$ and $H$ has a unique
orbital \mbox{$A= (a,b)$}.  
\be
\item If $H$ has a consistent controller, $c$,
then any element $g$ of $H$ which realizes both ends of $A$ actually
realizes $A$.
\item If $H$ has an inconsistent
controller $c$, then no element of $H$ realizes $A$.
\ee  
\ec 

Suppose that we know that $A$ is an orbital of a subgroup
$H\leq\ploi$.  Let $H_A$ be the set of elements of $\ploi$ each of
which is equal to the restriction of some element of $H$ on the
orbital $A$, and behaves as the identity off of $A$.  $H_A$ is
trivially a group with unique orbital $A$, and is a quotient of $H$.
We will call $H_A$ the \emph{projection of $H$ on $A$}.

We will now generalize our language somewhat.  Let $H$ be a subgroup
of $\ploi$ with an orbital $A$, and let $H_A$ be the projection of $H$
on $A$.  $H_A$ has a controller $\tilde{c}$ for $A$.  Let
$\rho:H\to H_A$ be the projection homomorphism on the orbital $A$.
Let 
\[T = \left\{c\in H\,|\,c\rho \textrm{ is a controller for $H_A$ on
$A$}\right\}.\]  We will call any element of $T$ \emph{a controller of
$H$ on $A$}\index{controller!on an orbital}.  Again, given an element
$c\in H$ which is a controller of $H$ on $A$, we can write elements of
$H$ in a unique $c$-form.  I.e., if $g\in H$, and $c$ is a controller
for $H$ on $A$, then there is an integer $k$ so that $g =
c^kg^{\!\!\!^{\circ}}$, where $g^{\!\!\!^{\circ}}$ will not realize
either end of $A$.

\section{Connecting the Two Algebraic Descriptions}
In this section, we will complete the proofs of our two main theorems
by introducing a new geometric technique.  Our technique allows us to embed
a solvable group of derived length $n$ in $\ploi$ into a more tractible
group in $\ploi$ with derived length $n$.
\subsection{Technical Preliminaries}
The next lemma is a technical lemma that we will use in completing our
proof of Theorem \ref{solveClassification}.

\bl
\label{oneBumpGenerators}

If $G$ is a solvable subgroup of $\ploi$ with derived length $n$,
generated by a collection $\Gamma$ of elements of $\ploi$ which each
admit exactly one orbital, and where no generator can be conjugated by
an element of $G$ to share an orbital with a different generator, then
$G$ is isomorphic to a group in the class $\left\{1\right\}\pc$ with
derived length $n$.

\el

pf: Before getting into the main body of the proof, note that the
hypothesies imply the main points of Lemma \ref{solvableQualities},
namely, that $G$ is balanced and admits no transition chains of length
two.  Further, the hypothesies further imply that each generator is
the only generator with that orbital, and that no element orbital is a
union of element orbitals that do not realize the original orbital.

We now enter the main body of the proof.  We will proceed by induction
on $n$.

If $n = 0$ then $G$ is the trivial group, and $G\in
\left\{1\right\}\pc$.  If $n = 1$ then $G$ is abelian, and in
particular, there can be at most countably many generators in
$\Gamma$, all of which have disjoint support, so that $G$ is
isomorphic with a countable (or finite) direct sum of $\Z$ factors, so
$G\in \{1\}\pc$.

Now let us suppose that $n > 1$ and that the statement of the lemma is
correct for any such solvable group with derived length $n-1$.  Let
$X$ represent the generators in $\Gamma$ whose orbitals are all depth
$2$ in $G$ or deeper.  Let $Y$ be the set of elements in $\Gamma$ whose
orbitals have depth $1$.  We note in passing that the cardinality of
$Y$ is at most countably infinite, and that the collection of orbitals
of the elements of $Y$ actually form the orbitals of the group $G$.
We will assume that all of the elements in $Y$ move points to the
right on their orbitals.  We can partition the elements in $X$ into
sets $P_y$, where the $y$ index runs over the elements in $Y$, and
where an element of $X$ is in $P_y$ if and only if that element's
orbital is contained inside the orbital of $y$.  Given $y\in Y$, if
$P_y$ is empty, then define $H_y = \langle y \rangle$.  Otherwise, let
$\gamma\in P_y$, and suppose that $\gamma$ has smallest depth possible
for the elements in $P_y$, and that $\gamma$ has orbital $A = (a,b)$.
$y$ has a fundamental domain $D_y = [a,ay)$, and each element of $P_y$
may be conjugated by some power of $y$ so that the resultant element's
orbital lies in the fundamental domain $D_y$ (if some element,
$\beta$, conjugates to contain $a$ in its orbital, then either that
conjugate has that its orbital fully contains the orbital $A$, which
is impossible by our choice of $\gamma$ as having a minimal depth
orbital of the orbitals of all the elements in $P_y$, or the signed
orbitals of $\gamma$ and of the conjugate of $\beta$ form a transition
chain of length two, which is impossible since $G$ is solvable).  We
can now replace $P_y$ by the conjugates of the original $P_y$ found
above so each element of $P_y$ has its orbital in $[a,a_y)$, and the
group generated by the new $P_y$ with $y$ will be identical to the
group generated by the old $P_y$ with $y$.  However, now that all of
the elements of $P_y$ have supports in the same fundamental domain of
$y$, we have that $H_y = \langle P_y, y\rangle$ is isomorphic to
$K_y\wr\Z$, where $K_y = \langle P_y \rangle$.  But $K_y$ is a
solvable group of precisely the type mentioned in the hypothesies of
the lemma, with derived length $k$ less than $n$, so that $K_y$ is
isomorphic to a group in $\left\{1\right\}\pc$ with derived length
$k<n$.  Therefore, $H_y$ is isomorphic to a group in
$\left\{1\right\}\pc$ (being the result of a group in
$\left\{1\right\}\pc$ being wreathed with a $\Z$ factor on the right)
with derived length $k+1\leq n$.  This argument holds for every $y$ in
$Y$, so $G\cong \bigoplus_{y\in Y}H_y$, where all of the groups in
this countable direct sum have derived length less than or equal to
$n$ (and at least one of them has derived length $n$).  In particular,
we see that $G$ is isomorphic to a group in $\left\{1\right\}\pc$ with
derived length $n$.  \qquad$\diamond$

\subsection{The Split Group}
One new tool for technical analysis of a subgroup $G$ of $\ploi$ is
the split group of $G$.  It is motivated by the hypothesies of Lemma
\ref{oneBumpGenerators}.  We define the split group of a subgroup of
$\ploi$ below.

Let $P_s$ represent the set of subgroups of $\ploi$.  We define a
function $S:P_s\to P_s$.  Given a group $G$ which is a subgroup of
$\ploi$, fix the notation $\Gamma_G$ to represent the maximal set of
one-orbital elements of $\ploi$ such that if $\gamma\in\Gamma_G$, then
$\gamma$ is identical to an element $g\in G$ over $\gamma$'s orbital.
Define $S$ by the rule that $S(G)= \langle \Gamma_G\rangle$ (if $G$ is
the trivial group, define $S(G) =G$).  Given a group $G\in P_s$, we
will call the group $S(G)$ the \emph{split group of $G$}.  Note that
$G\leq S(G)$.  In this section, we will analyze further properties of
the split group $S(G)$, in the case that $G$ is solvable.

Throughout the remainder of the section, given a group $G$, and any
$\ga\in\Gamma_G$, let $A_{\ga}$ always denote the orbital of $\ga$.
Below, we mention some further basic facts relevant to the analysis of
the split group of a group, some using the notation just established.
We leave the proofs to the reader.

\brk
\label{splitProperties}

Suppose $G$ is a solvable subgroup of $\ploi$, then
\be

\item The depth of any signed orbital of $G$ is a positive integer.
\item If $A=(a,b)\subset I$, and both $g$, $h\in G$ have orbital $A$,
then the depths of the signed orbitals $(A,\,g)$ and $(A,\,h)$ of $G$
are the same.
\item If $(A,g)$ is a signed orbital of $G$ with depth $k$, and $h\in G$, then
$(Ah,g^h)$ is a signed orbital of $G$ with depth $k$.  
\item If $\ga$, $\gb\in\Gamma_G$ and the orbital of $\ga$ is contained in
the orbital of $\gb$, then $\ga^{\gb}$ is in $\Gamma_G$. 
\item If $\ga$, $\gb\in\Gamma_G$, then
$\left\{(A_{\ga},\,\ga),\,(A_{\gb},\,\gb)\right\}$ is not a transition
chain of length two for $S(G)$.

\ee 

\erk

Note that in the case where $G$ is solvable, the first and second
points above allow us to define a \emph{$G$-depth} for elements in
$\Gamma_G$.  Namely, if $\gamma\in\Gamma_G$ and the orbital of
$\gamma$ is $A$, then there is a positive integer $k$ so that all
signed orbitals of the form $(A,\,h)$ for $G$ (in particular, this
includes all of the signed orbitals where the element $h$ agrees with
$\gamma$ over $A$) have depth $k$, so we can say $\gamma$ has
$G$-depth $k$.

The following lemma follows easily from the definition above, and
points $4$ and $5$ of the previous remark. It provides us with
something like a normal form for writing elements of $S(G)$ in terms
of the generators in $\Gamma_G$.

\btl 
\label{splitForm}

Suppose $G$ is a solvable subgroup of $\ploi$ and $h\in S(G)$, then
given a product decomposition $h=\tau_1\tau_2\cdots\tau_k$ so that
each $\tau_i$ is an element of $\Gamma_G$ for each positive integer
$\,i\leq k$, then there is a product decompostion $h =
\theta_1\theta_2\cdots\theta_k$ with each $\theta_i\in\Gamma_G$, so
that whenever $r$ and $s$ are positive integers so that $r<s \leq k$,
we have the $G$-depth of $\theta_s$ is greater than or equal to the
$G$-depth of $\theta_r$.

\etl 

pf: 

We will show that if $h\in S(G)$, and $h = \alpha_1\alpha_2$ with
$\ga_1$, $\ga_2\in\Gamma_G$, then we can write $h = \theta_1\theta_2$
with $\theta_1$ and $\theta_2\in\Gamma_G$ and with the $G$-depth of
$\theta_2$ greater than or equal to the $G$-depth of $\theta_1$, and
with the set of $G$-depths of the $\theta_i$'s the same as the set of
$G$-depths of the $\ga_i$'s.  We then will have our lemma by re-writing
longer products, improving the product locally on adjacent pairs of
length two until there are no more improvements to be made.

Let us assume the $G$-depth of $\ga_1$ is greater than the $G$-depth
of $\ga_2$, as the other cases are trivial.  In this case let
$\theta_1 = \ga_2$ and let $\theta_2 = \ga_1^{\ga_2}$.  If $\ga_2$ and
$\ga_1$ have orbitals that intersect, then point $5$ from Remark
\ref{splitProperties} and the definition of $G$-depth guarantee that
the orbital of $\ga_1$ is contained in the orbital of $\ga_2$, and
then point $4$ of Remark \ref{splitProperties} gives us that
$\theta_2\in\Gamma_G$.  In the case that the orbitals of $\ga_1$ and
$\ga_2$ are disjoint, then $\theta_2 = \ga_1\in\Gamma_G$, so in all
cases $\theta_1$ and $\theta_2\in\Gamma_G$.  Now we can afford the
cost of pushing $\ga_2$ past $\ga_1$:

\[
h =\ga_1\ga_2 = \ga_2\cdot\ga_2^{-1}\ga_1\ga_2 = \theta_1\cdot\theta_2
\]

Note that the $G$-depth of $\theta_1$ is the $G$-depth of $\ga_2$ and the
$G$-depth of $\theta_2$ is the $G$-depth of $\ga_1$ (by a variant of
point three from Remark \ref{splitProperties}).

\qquad$\diamond$

We will need one more technical lemma before we can come to some
meaningful conclusions about the split group of a solvable group $G$.

\btl \label{firstOrbitalDominance} 

Suppose $G$ is a solvable subgroup of
$\ploi$.  If $\left\{\ga_1,\,\ga_2,\,\ldots,\,\ga_k\right\}$ is a
collection of elements of $\Gamma_G$ for some integer $k$, and for
each integer $j$ with $2\leq j \leq k$ we have that
$\overline{A}_{\ga_j}\subset A_{\ga_1}$, then the product $\tau_{1,k}
= \ga_1\ga_2\cdots\ga_k$ will have support $A_{\ga_1}$.

\etl

pf:

Call a finite product of elements of $\Gamma_G$ a \emph{first orbital
dominant product} if the closure of each orbital of the later terms in
the product is contained in the orbital of the first term of the
product.

It is immediate that the support of a first orbital dominant product
is contained in the orbital of the first element of the product.  We
will now show that such a product admits no fixed points in the
orbital of the first element of the product.

Let us suppose first that $k$ is the smallest integer so that a first
orbital dominant product of $k$ elements of $\Gamma_G$ has a fixed
point in the orbital of the first element of the product, and let
$\ga_1\ga_2\cdots\ga_k$ be such a product.  We will derive a
contradiction, and therefore we will be able to conclude that any
first orbital dominant product will have the whole of the orbital of
the first element in the product as its support.

We note in passing that $k>1$.

By reference to Technical Lemma \ref{splitForm}, we assume (without
risking the first orbital dominant nature of the product) that if
$i<j$ are integers with $1\leq i <j\leq k$ then the $G$-depth of
$\ga_i$ is less than or equal to the $G$-depth of $\ga_j$.  Let us
build some notation for subproducts of our newly re-arranged
$\ga_i$'s.  Given any integers $i$ and $j$ so that $1\leq i\leq j\leq
k$ let $\tau_{i,j}$ represent the product
$\ga_i\ga_{i+1}\cdots\ga_{j}$.

If $x\in A_{\ga_1}$ so that $x\tau_{1,\,k} = x$, then since $x$ is not
a fixed point of $\tau_{1,\,k-1}$ we see that $x\in A_{\ga_k}$.  Since
the fixed set of a product of elements of $\ploi$ is closed, let us
further assume that $x$ is the infimum of the points in $A_{\ga_1}$
that are fixed by $\tau_{1,\,k}$.  Denote $x$ by $y_0$, and further,
let $y_j=y_0\tau_{1,\,j}$ for each integer $j$ with $1\leq j\leq k$,
so that, for example, $y_k = x$.  If there is a positive integer
$j\leq k$ with $y_j = y_{j-1}$ then the product
$\ga_1\ga_2\cdots\ga_{j-1}\ga_{j+1}\ga_{j+2}\cdots\ga_k$ will also
have fixed point $x$, contradicting the minimality of $k$, so we see
that for each positive integer $j$ with $j \leq k$, $y_{j-1}\neq y_j$.
It is immediate from the last condition that for any $j$ with $0\leq
j<k$, $y_j\in A_{\ga_{j+1}}$.  The previous sentence implies that for
each $i$ with $1\leq i<k$ we have $A_{\ga_{i+1}}\cap A_{\ga_{i}}\neq
\emptyset$.  Therefore, $A_{\ga_{i+1}}\subset A_{\ga_{i}}$ since the
$G$-depths of the $\ga_i$ are non-decreasing and $G$ admits no
transition chains of length two.

By the definition of $\Gamma_G$, for each integer $i$ with $1\leq
i\leq k$ there are elements $g_i\in G$ with $\ga_i$ behaving like
$g_i$ over $A_{\ga_i}$.  We observe that the product $g = g_1g_2\cdots
g_k$ will also fix $x$, since $y_{i-1}\in A_{\ga_{i}}$ for all integer
$i$ with $1\leq i\leq k$.  Denote by $a_k$ the left hand endpoint of
$A_k$.  We will show that $a_kg\neq a_k$, implying that $g$ has an
orbital containing the left end of the orbital $A_{\ga_k}$ of $g_k$.

We observe immediately that $(a_k,x)\subset A_{\ga_{k-1}}$.  Define
$b_k = a_k$ and inductively define, for each integer $j$
with $0\leq j<k$ the sets $(b_j,y_j) =
(b_{j+1},y_{j+1})\ga_{j+1}^{-1}$. Since $(b_j, y_j)\subset
A_{\alpha_{j+1}}$ for all $j\in\N$ with $j < k$, we
see that $g$ agrees with $\tau_{1,k}$ over $(b_0,y_0)$, and therefore
by the continuity of $g$ and $\tau_{1,k}$, these two maps agree over
$[b_0,y_0]$.  But $y_0 = x$, where $g$ is fixed, while $b_0\tau_{1,k} =
a_k$.  Since $\tau_{1,k}$ is not fixed at $b_0$ (this follows from the
facts that $k>1$ and $\overline{A}_j\subset A_1$, for $j\in \N$ with
$2\leq j \leq k$), $g$ moves $a_k$. In
particular, $g$ has an orbital $B$ containing $a_k$ but not $x$, so
that $\left\{(B,\,g),\,(A_k,\,g_k)\right\}$ is a transition chain of
length two for $G$, which contradicts the solvability of $G$.

\qquad$\diamond$

We are now in a position to analyze the depth of the split group of a
solvable group $G$.  \bl
\label{splitStable}

Suppose $G$ is a solvable subgroup of $\ploi$ that has derived length
$u$ for some $u\in\N$.  If $(A,h)$ is a signed orbital of $S(G)$, then
there is $g\in G$ so that $(A,g)$ is a signed orbital of $G$.

\el
pf:

If $u=0$ there is nothing to prove, so let us assume that $u>0$, so
that $S(G)$ is not the trivial group.  In particular, we assume that
$S(G)$ has some associated signed orbitals.

Suppose $(A,h)$ is a signed orbital of $S(G)$, and let $h =
\ga_1\ga_2\cdot\cdots\cdot\ga_k$ for some positive integer $k$, where
each $\ga_i$ is in $\Gamma_G$.  In this argument, we will give an
algorithm that steadily improves $h$ and the product expression for
$h$ (always preserving the fact that $h$ has orbital $A$).  After
each improvement, we will assume the product is re-indexed as in the
product above, so that, for instance, $k$ will always refer to the
length of the current product.  At some point, we will see that $G$
has an element with orbital $A$.

Our algorithm begins in the next paragraph.  At some points we improve
the product definition of $h$, and we then begin the algorithm again.
This could easily lead to infinite loops in our algorithm.  After the
algorithm is fully specified, we will argue that the direction to
restart the algorithm, given in the algorithm's final paragraph,
avoids the creation of such an infinite loop.  We will leave the
argument for termination of the algorithm for each of the previous
directions to re-start to the reader.

\vspace{.1 in}

\underline{Algorithm Initiation}
\be
\item Assurance that $k$ is large.

If $k = 1$ then $A_{\ga_1} = A$ and $G$ has an element $g$ which
behaves as $\ga_1$ on $A$.  Therefore let us assume that $k>1$.

\item $G$-depth ordering, and some notation.

By Technical Lemma \ref{splitForm} we can assume that the product
decomposition for $h$ has the property that if $r$ and $s$ are
positive integers so that $r<s\leq k$, then the $G$-depth of the
orbital of $\ga_s$ is greater than or equal to the $G$-depth of $\ga_r$.

Given integers $i$ and $j$ with $1\leq i\leq j\leq k$ set $\tau_{i,j}
= \ga_i\ga_{i+1}\cdots\ga_j$, so that $h = \tau_{1,k}$.

\item $A\subset A_{\ga_1}$.

Suppose $A_{\ga_1}$ is disjoint from $A$.  In this case $\tau_{2,k}$
must have orbital $A$. Replace $h$ by $\tau_{2,k}$, and start the
algorithm again using the one-shorter product description of the new
$h$, re-indexed.  We can assume from here out that $A\cap
A_{\ga_1}\neq\emptyset$.  We will refer to $A_{\ga_1}$, as determined
at this stage, to be the \emph{leading orbital of the product}.

It is immediate from the definition of $G$-depth that if $r$, $s$ are
integers so that $1\leq r<s \leq k$ then either $A_{\ga_s}\subset
A_{\ga_r}$ or these two element orbitals are disjoint (recall point
$5$ of Remark \ref{splitProperties}).

We can strengthen the result of the paragraph two paragraphs back by
noting that if $A$ is not contained in $A_{\ga_1}$, then some end $e$
of $A_{\ga_1}$ must be in $A$.  In this case $e$ is moved by one of
the $\ga_i$ for an integer $i$ with $2 \leq i\leq k$.  This last is
impossible by the previous paragraph, therefore $A\subset
A_{\ga_1}$.

\item Forcing $A_{\ga_i}\subset A_{\ga_1}$ for $1\leq i \leq k$.

Now, for any integer $i$ with $2\leq i\leq k$, the orbital $A_{\ga_i}$
must either be contained in $A_{\ga_1}$ or must be disjoint from
$A_{\ga_1}$.  In particular, if we set $C_j = A\tau_{1,j}$ for each
integer $j$ with $1\leq j \leq k$, we see that for each such $j$,
$C_j\subset A_{\ga_1}$.  But this implies that we can drop any $\ga_j$
from the product producing $h$, if that $\ga_j$ has support disjoint
from $A_{\ga_1}$, and the resulting product will still act as $h$ over
$A$.  If there are any such $j$, let us drop the corresponding $\ga_j$
from the product, re-index, and begin the algorithm again (to catch
the case that the resulting product has length one, for instance).  If
there are no such $j$, then we see that for each integer $j$ with
$2\leq j\leq k$ we must have $A_{\ga_j} \subset A_{\ga_1}$.

\item Finding how many of the $A_{\ga_i}$ are equal to $A_{\ga_1}$.

For each integer $i$ with $1\leq i\leq k$ note that there is an
element $g_i\in G$ so that $g_i$ behaves as $\ga_i$ over the interval
$A_{\ga_i}$.  

Suppose that for some integer $i$ in $2\leq i\leq k$ we have that
$A_{\ga_i}$ shares an end with $A_{\ga_1}$.  Since $G$ is balanced,
either $A_{\ga_i} = A_{\ga_1}$ or $g_i$ has at least two orbitals
contained in $A_{\ga_1}$, one sharing one end of $A_{\ga_1}$, and the
other sharing the other end of $A_{\ga_1}$.  It is easy in the second
case to show that $G$ admits transition chains of length two.  In
particular, we can now assume that for any integer $i$ where $2\leq i
\leq k$, the orbital $A_{\ga_i}$ either equals $A_{\ga_1}$ or has
closure contained in $A_{\ga_1}$.  

By the fact that the $G$-depth of the elements $\ga_i$ is
non-decreasing in $i$, we see that there is an integer $s$ with $1\leq
s\leq k$ so that for all integers $1\leq r\leq s$ we have $A_{\ga_r} =
A_{\ga_1}$ while for all integers $t$ with $s< t\leq k$ we see that
$\overline{A}_{\ga_t}\subset A_{\ga_1}$.

If $s = 1$, Technical Lemma \ref{firstOrbitalDominance} assures us
that $A = A_1$.  In this case, $g_1$ has orbital $A$ and we are
finished.  Therefore, let us assume that $s>1$, so that $A_{\ga_1} =
A_{\ga_2}$.

\item Collapsing the product $\ga_1\ga_2$ over $A_{\ga_1}$.

If $\tau_{1,2}$ is the identity, then $\tau_{3,k}$ behaves as $h$ over
$A$.  In this case, remove $\ga_1$ and $\ga_2$ from the product,
re-index, and begin the algorithm again.  If $\tau_{1,2}$ is not the
identity, continue below.

Let $\left\{B_1,\, B_2,\,\ldots,\,B_m\right\}$ be the orbitals of
$\tau_{1,2}$.  These orbitals are all contained in
$A_{\ga_1}$.  Note that $\tau_{1,2}$ behaves as an element of $G$ over
$A_1$.  In particular, there are elements $\gb_i\in\Gamma_G$ for each
integer $1\leq i\leq m$ so that $\gb_i$ behaves as $\tau_{1,2}$ over
$B_i$.  

Let $\Upsilon = \left\{i\in\N\,|\,1\leq i \leq m, B_i\cap
A\neq\emptyset\right\}$.  Form the product
$\omega = (\prod_{i\in\Upsilon}\gb_i)\ga_3\ga_4\cdots\ga_k$, which is a product
of elements of $\Gamma_G$ which behaves as $h$ over $A$, but has fewer
elements in the product with orbital $A_{\ga_1}$.  Begin the algorithm again
using the product definition of $\omega$ to define the new $\ga_i$,
indexed appropriately.

\ee

\underline{Algorithm Termination}

\vspace{.1 in}

We will now argue that the algorithm given above terminates, in terms
of the final paragraph's re-direction to repeat the algorithm.  The
previous re-directions cannot create an infinite loop for
straightforward reasons which we leave to the reader.

The final paragraph of the algorithm either decreases the number of
elements in the product which have orbital equal to the leading
orbital of the product (which can only happen a finite number of times
for any particular leading orbital, bounded by $s-1$ for the initial
$s$ of that leading orbital), or it will decrease the size of the
leading orbital.  If the algorithm decreases the size of the leading
orbital $u$ times then $G$ will admit a tower of height $u+1$, which
is impossible as $G$ has derived length $u$.  

The previous paragraph implies that the algorithm must terminate
before the steps in the final paragraph of the algorithm description
decrease the size of the leading orbital for the $u$'th time, which
means that it will have found an element of $G$ with orbital $A$.

\qquad$\diamond$

\bc 
\label{splitLength}

Suppose $G$ is a subgroup of $\ploi$.  The derived length of $G$
equals the derived length of $S(G)$.  

\ec

pf: 

If $G$ is solvable then Lemma \ref{splitStable} applies; given any
tower $T$ of $S(G)$ we can find a tower with the same orbitals for
$G$.  Hence, the derived length of $S(G)$ is less than or equal to the
derived length of $G$.  But since $G\leq S(G)$ we see that the derived
length of $G$ is automatically less than or equal to the derived
length of $S(G)$.  In particular, if $G$ is solvable, these derived
lengths are equal.

If $G$ is non-solvable, then both $G$ and $S(G)$ have towers of
arbitrary height, so that both groups are unsolvable.  

\qquad$\diamond$

\subsection{Split Groups and $\mm$}

The following lemma, and its corollary, complete our proof of both
Lemma \ref{shallowGroupEmbeddings} and Theorem
\ref{solveClassification}.  Note that the corollary is simply a
restatement of Lemma \ref{solveInM}.

\bl 
\label{HInOneP}

If $G$ is a solvable subgroup of $\ploi$, then $S(G)$ is isomorphic to
a group in $\left\{1\right\}\pc$.

\el 

pf: 

Suppose that $n\in\N$ and that $G$ is a solvable subgroup of $\ploi$
with derived length $n$.  By Corollary \ref{splitLength} we know that the
derived length of $S(G)$ is also $n$.  By Lemma \ref{solvableQualities}
we see that $G$ admits no transition chains of length $2$.

Let $X_1$ represent the set of signed orbitals of $H$ with depth $1$.
For each orbital $A$ in $O_{X_1}$, there is a non-empty set of
controllers of $H$ for $A$ in the set $\Gamma$.  Let
$\phi_1:O_{X_1}\to \Gamma$ represent a function that associates to
each orbital in $O_{X_1}$ a controller of $H$ for that orbital which
moves points to the right on the orbital.  Let $Y_1 = O_{X_1}\phi_1$
be the image of $\phi_1$.  We note that each pair of elements in $Y_1$
have disjoint support, and trivially, that no element of
$Y_1$ can be conjugated by an element of $H$ to share an orbital with
a different element of $Y_1$.  Now, by the definition of controller,
$Y_1$ consists of a set of generators in $\Gamma$ sufficient so that
the set $\Gamma_1 = Y_1 \cup (\Gamma \backslash S_{X_1})$ generates
$H$.

For each $y\in Y_1$, let $A_1^y$ represent the orbital of $y$.  We may
partition the elements of $\Gamma_1$ into sets $P^y_1$ indexed by the
set $Y_1$ so that $\gamma\in P^y_1$ if the orbital of $\gamma$ is
contained in $A_1^y$.  Now let $X^y_2$ represent the set of signed
orbitals of elements in $P^y_1$ with depth two in $H$.  If $X^y_2$ is
not empty, let $(A^y_2,\gamma^y_2)$ be an element in $X^y_2$, and let
$a_y$ be the left end of the orbital $A^y_2$.  Since the orbitals of
elements of $H$ are always orbitals of $G$ by Lemma \ref{splitStable},
we see that $H$ admits no transition chains of length two.  In
particular, we can use $y$ to conjugate every signed orbital of
$X^y_2$ into the fundamental domain $[a_y,a_yy)$ to produce the set
$D^y_2$ of signed orbitals.  Likewise, we can conjugate every signed
orbital in $P^y_1$ of depth greater than two into the fundamental
domain as well, producing the set $T^y_2$. The collections of
conjugates in the fundamental domain have nice properties:

\be

\item If two elements of $D^y_2$ have orbitals that non-trivially
intersect each other, then they actually have identical orbitals.
\item $H_{A^y_1}=\langle y,S_{D^y_2},S_{T^y_2}\rangle$.  

\ee 

Let $\phi^y_2:O_{D^y_2}\to S_{D^y_2}$ be a function that picks for
each orbital of depth two in $O_{D^y_2}$ a controller that moves
points to the right for that orbital as before.  Let $Y_2$ be the
union of all the images of the functions $\phi^y_2$ across the $Y$
index set, so that we have now picked a controller that moves points
to the right on its orbital for every conjugacy class of depth two
orbitals of $H$ (note that a conjugate of a controller is also a
controller), so that if $\Gamma_2 = (\Gamma_1\backslash (\cup_{y\in
Y}S_{X^y_2}))\cup Y_2$ then $H = \langle \Gamma_2\rangle$.

In a like fashion we can inductively proceed to pick sets of
controllers, one for each conjugacy class of element orbital of depth
$i$, where $i$ is an index less than or equal to $n$, in exactly the
same fashion as discussed for forming the set $Y_2$ above.  This
process will steadily improve the sets of generators $\Gamma_i$, so
that finally $\Gamma_n$ will be a set of generators for $H$ where each
generator in $\Gamma_n$ has exactly one orbital, and where no
generator can be conjugated in $G$ to share an orbital with another
generator in $\Gamma_n$.  In particular, by Lemma
\ref{oneBumpGenerators}, $H$ is isomorphic to a group in
$\left\{1\right\}\pc$ with derived length $n$, and therefore $G$ is
isomorphic to a subgroup of a group in the class $\left\{1\right\}\pc$
with derived length $n$.  \qquad$\diamond$

\bc
\label{solveEmbeddings}
If $G$ is solvable in $\ploi$ of derived length $n$, then $G$ embeds
in $G_n$.  
\ec 

pf: 

$G$ embeds in $S(G)$, which also has derived length $n$.  $S(G)$
is isomorphic to a group in $\left\{1\right\}\pc$ with derived length
$n$ by Lemma \ref{HInOneP}.   Now, by Lemma \ref{BWinM}, $S(G)$
embeds in $G_n$.  \qquad$\diamond$

\bibliographystyle{amsplain}
\bibliography{ploiBib}
\end{document}